\documentclass[11pt]{article}

\usepackage{amsmath,amssymb,amsthm,mathtools,mathrsfs,braket,bbm,dsfont}
\usepackage{enumitem}

\usepackage{authblk}
\usepackage{booktabs}

\usepackage[numbers,sort&compress]{natbib}
\usepackage{multirow}
\usepackage[hidelinks]{hyperref}
\usepackage{textgreek}
\usepackage[T1]{fontenc}

\usepackage{interval}
\intervalconfig{soft open fences}

\newtheorem{theorem}{Theorem}
\newtheorem{lemma}[theorem]{Lemma}

\newtheorem{remark}[theorem]{Remark}
\newtheorem{corollary}[theorem]{Corollary}
\newtheorem{definition}[theorem]{Definition}

\newcommand{\acks}[1]{\section*{Acknowledgements}#1}
\usepackage{changes}

\title{Linear Convergence in Hilbert's Projective Metric for Computing Augustin Information\\and a R\'{e}nyi Information Measure}

\author[1]{Chung-En Tsai$^\ast$}
\author[2]{Guan-Ren Wang$^\ast$}
\author[3,4,5,6,7]{\authorcr Hao-Chung Cheng}
\author[1,2,4,5]{Yen-Huan Li}

\affil[1]{Department of Computer Science and Information Engineering,\protect\\National Taiwan University}
\affil[2]{Graduate Institute of Networking and Multimedia, \protect\\ National Taiwan University}
\affil[3]{Department of Electrical Engineering and Graduate Institute of\protect\\Communication Engineering, National Taiwan University}
\affil[4]{Department of Mathematics, National Taiwan University}
\affil[5]{Center for Quantum Science and Engineering, \protect\\ National Taiwan University}
\affil[6]{Physics Division, National Center for Theoretical Sciences}
\affil[7]{Hon Hai (Foxconn) Quantum Computing Centre}

\date{\vspace{-5ex}}

\newcommand{\R}{\mathbb{R}}
\newcommand{\N}{\mathbb{N}}
\newcommand{\E}{\mathbb{E}}

\DeclareMathOperator{\ri}{ri}

\newcommand{\abs}[1]{\left\lvert#1\right\rvert}
\newcommand{\norm}[1]{\left\lVert#1\right\rVert}
\newcommand{\dH}{d_{\mathrm{H}}}

\newcommand\blfootnote[1]{%
  \begin{NoHyper}%
  \renewcommand\thefootnote{}\footnote{#1}%
  \addtocounter{footnote}{-1}%
  \end{NoHyper}%
}

\begin{document}

\maketitle

\begin{abstract}

Consider the problems of computing the Augustin information and a R\'{e}nyi information measure of statistical independence, previously explored by Lapidoth and Pfister (\textit{IEEE Information Theory Workshop}, 2018) and Tomamichel and Hayashi (\textit{IEEE Trans. Inf. Theory}, 64(2):1064–-1082, 2018). 
Both quantities are defined as solutions to optimization problems and lack closed-form expressions.
This paper analyzes two iterative algorithms: Augustin’s fixed-point iteration for computing the Augustin information, and the algorithm by Kamatsuka et al. (arXiv:2404.10950) for the R\'{e}nyi information measure.
Previously, it was only known that these algorithms converge asymptotically. 
We establish the linear convergence of Augustin’s algorithm for the Augustin information of order $\alpha \in (1/2, 1) \cup (1, 3/2)$ and Kamatsuka et al.’s algorithm for the R\'{e}nyi information measure of order $\alpha \in [1/2, 1) \cup (1, \infty)$, using Hilbert’s projective metric.
\\\\
\end{abstract}

\blfootnote{\hspace{-15pt}$^\ast$Both authors contribute equally to this work.}


\section{Introduction}\label{sec:introduction}

Denote by $\Delta ( [d] )$ the set of probability distributions over the finite set $[d] \coloneqq \set{ 1, \ldots, d }$. 
For any $\alpha \in \interval[open right]{0}{1} \cup \interval[open]{1}{\infty}$, the order-$\alpha$ Augustin information is defined by the following optimization problem \citep{Augustin1978}:
\begin{equation}\label{eq:augustin_information}
	\min_{x\in\Delta([d])} f_{\mathrm{Aug}}(x),
	\quad f_{\mathrm{Aug}}(x):= \E_{p\sim P}\left[ D_\alpha( p \parallel x ) \right],
\end{equation} 
where $P$ is a given probability distribution over $\Delta([d])$, and
\[
	D_\alpha(p \parallel q) \coloneqq \frac{1}{\alpha - 1}\log\sum_{s\in S} p(s)^{\alpha} q(s)^{1-\alpha},
	\quad\forall p,q\in\Delta(S)
\]
is the order-$\alpha$ R\'{e}nyi divergence.
The Augustin information characterizes, e.g., the cutoff rate, the strong converse exponent, and the error exponent in the channel coding problem \citep{Arimoto1973,Csiszar1995,Csiszar2011,Nakiboglu2019,Wang2024}.
When $\alpha = 0$, the optimization problem \eqref{eq:augustin_information} becomes a special case of the log-optimal portfolio \citep{Cover1984}, which is equivalent to the definition of the maximum-likelihood estimate in Poisson inverse problems 
\citep{Vardi1993}.


The optimization problem \eqref{eq:augustin_information} does not admit a closed-form expression. 
While the optimization problem is convex, the objective function violates the standard smoothness assumption in the optimization literature. 
Therefore, even the convergence guarantees of projected gradient descent, arguably the simplest convex optimization algorithm, do not directly apply \citep{You2022}. 

\citet{Augustin1978} proposed the 
following 
fixed-point iteration to solve the optimization problem \eqref{eq:augustin_information}: 
\begin{equation}\label{eq:cover}
	x_{t+1} = x_t \odot (-\nabla f_{\text{Aug}}(x_t)),
	\quad \forall t \in \mathbb{N} .  
\end{equation}
where 
$\odot$ denotes the entry-wise product.
The algorithm was later rediscovered by \citet{Karakos2008}\footnotemark{}. 
When $\alpha = 0$, this fixed-point iteration coincides with Cover's method for computing the log-optimal portfolio \citep{Cover1984}, and is equivalent to the expectation maximization algorithm for solving Poisson inverse problems \citep{Richardson1972,Lucy1974,Shepp1982,Vardi1993}.

\footnotetext{\citet{Karakos2008} proposed an alternating minimization method whose iteration consists of two steps.
Combining the two steps yields Augustin's fixed-point iteration.}

Recently, \citet{Kamatsuka2024} proposed an algorithm similar to Augustin's fixed-point iteration to compute a R\'{e}nyi information measure 
of statistical independence, which was explored 
by \citet{Lapidoth2019} and \citet{Tomamichel2018}.
For 
any 
$\alpha\in[0,1)\cup(1,\infty)$, this order-$\alpha$ R\'{e}nyi information measure is defined 
by the following optimization problem: 
\begin{equation}\label{eq:renyi_information_measure}
	\min_{x\in\Delta([m])}\min_{y\in\Delta([n])} f_{\mathrm{Ren}}(x,y),
	\quad f_{\mathrm{Ren}}(x,y) 
	\coloneqq
	D_\alpha( p \parallel x \otimes y ),
\end{equation}
where $p$ is a given probability distribution over $[m] \times [n]$ and $\otimes$ 
denotes the tensor product. 
The R\'enyi information measure emerges in the error exponent of a hypothesis testing problem, where we test against the independence of two random variables given 
independent and identically distributed (i.i.d.)
samples from their joint distribution \citep{Lapidoth2018,Lapidoth2019,Tomamichel2018}.

\citeauthor{Kamatsuka2024}'s algorithm to compute the R\'{e}nyi information measure iterates as: 
\begin{equation}\label{eq:alternating_minimization}
\begin{alignedat}{2}
	x_{t+1} &= Z_{1, t}^{-1} \cdot x_t \odot ( -\nabla_x f_{\mathrm{Ren}}(x_t, y_t))^{1/\alpha}, \\
	y_{t+1} &= Z_{2, t}^{-1} \cdot y_t \odot ( -\nabla_y f_{\mathrm{Ren}}(x_{t+1}, y_t))^{1/\alpha},
\end{alignedat}
\end{equation}
where $Z_{1, t}$ and $Z_{2, t}$ are normalizing constants, 
ensuring that $x_{t+1}$ and $y_{t+1}$ remain probability distributions. 
The notation $v^r$ denotes the entry-wise power for any vector $v$ and number $r$. 
This iterative algorithm is reminiscent of Augustin’s fixed-point iteration but differs in the powers applied to the gradients.

The convergence behaviors of Augustin's fixed point iteration and Kamatsuka et al.'s algorithm remain largely unclear. 
For Augustin's fixed-point iteration, \citet{Karakos2008} and \citet{Nakiboglu2019} have shown that it asymptotically converges for $\alpha \in \interval[open]{0}{1}$; 
\citet{Iusem1992} and \citet{Lin2021} have proved a convergence rate of $O ( 1 / t )$ for the case where $\alpha$ approaches zero. 
For Kumatsuka et al.'s algorithm, \citet{Kamatsuka2024} have shown that it asymptotically converges for $\alpha \in \interval[open right]{1 / 2}{1} \cup \interval[open]{1}{\infty}$. 

We aim to carry out non-asymptotic analyses for the two algorithms.
One common approach to analyzing an iterative method is to show that it is 
contractive
under a suitable metric.
Since the two algorithms \eqref{eq:cover} and \eqref{eq:alternating_minimization} map positive vectors to positive vectors, we view them as positive dynamical systems and consider the so-called Hilbert's projective metric \citep{Lemmens2012,Krause2015}.


In this work, we prove that with respect to Hilbert's projective metric, Augustin's fixed-point iteration is contractive for $\alpha\in(1/2,1) \cup (1,3/2)$, and 
\citeauthor{Kamatsuka2024}'s algorithm is also contractive for $\alpha\in(1/2,1) \cup (1,\infty)$.
Based on these contractivity results, we establish the following non-asymptotic convergence guarantees for the two algorithms.

\begin{itemize}
\item For computing the Augustin information of order $\alpha\in \interval[open]{1/2}{1} \cup \interval[open]{1}{3/2}$, Augustin's fixed-point iteration converges at a rate of $O((2\abs{1-\alpha})^t)$ with respect to Hilbert's projective metric.
This improves on the previous asymptotic convergence guarantee \citep{Karakos2008,Nakiboglu2019} when $\alpha\in(1/2,1)$ and extends the range of convergence to include $\alpha\in(1,3/2)$.

\item For computing the R\'{e}nyi information measure of order $\alpha\in(1/2,1)\cup(1,\infty)$, the iterative algorithm of \citeauthor{Kamatsuka2024} converges at a rate of $O(\abs{1-1/\alpha}^{2t})$ with respect to Hilbert's projective metric.
When $\alpha=1/2$, this method also converges linearly if $p$ has full support.
This improves on the previous asymptotic convergence guarantee \citep{Kamatsuka2024}.
\end{itemize}

\paragraph{Notations}
We write $\mathbb{R}_+$ and $\mathbb{R}_{++}$ for the sets of non-negative and strictly positive numbers, respectively. 
For any positive integer $n$, we write $[n]$ for the set $\set{ 1, \ldots, n }$. 
Let $v \in \mathbb{R}^d$ and $A, B \in \mathbb{R}^{m \times n}$. 
We write $v ( i )$ for the $i$-th entry of the vector $v$, and $A ( i, j )$ the $(i, j)$-th entry of the matrix $A$. 
We write $A \odot B$ for the entry-wise product between $A$ and $B$. 
We write $A^r$ for the matrix $( A(i, j)^r )_{1 \leq i \leq m, 1 \leq j \leq n}$. 
For a set $S\subseteq\R^d$, we denote by $\ri S$ its relative interior.
We will adopt the convention that $0^0=0$, $0/0=1$, $\infty\cdot\infty=\infty$, $a\cdot\infty=\infty$ for any $a>0$, and $\log\infty=\infty$.
We call $\Delta([d])$ the probability simplex and view elements in $\Delta([d])$ as $d$-dimensional vectors.


\section{Related Work} \label{sec:related_work}

We have discussed Augustin's fixed-point iteration and \citeauthor{Kamatsuka2024}'s algorithm in Section~\ref{sec:introduction}.
This section reviews other optimization algorithms for computing the Augustin information and the R\'{e}nyi information measure.

For computing the Augusitin information of order $\alpha$, entropic mirror descent with Armijo line search \citep{Li2018} and with the Polyak step size \citep{You2022}, as well as a variant of Augustin's fixed-point iteration explored by \citet[Lemma 6]{Cheng2021}, all achieve asymptotic convergence for all $\alpha \in \interval[open right]{0}{1} \cup \interval[open]{1}{\infty}$. 
Riemannian gradient descent with the Poincar\'e metric \citep{Wang2024} converges at a rate of $O ( 1 / t )$ for all $\alpha \in \interval[open right]{0}{1} \cup \interval[open]{1}{\infty}$. 
An alternating minimization method due to \citet{Kamatsuka2024}\footnote{\citet{Kamatsuka2024} only claimed an asymptotic convergence guarantee in their paper.
We find that their Lemma~2 indeed implies a convergence rate of $O(1/t)$.} also achieves a converges rate of $O ( 1 / t )$, but for a narrower range of $\alpha \in \interval[open]{1}{\infty}$. 
None of the existing works have yet established a linear convergence rate. 

 

For computing the R\'enyi information measure of order $\alpha$, entropic mirror descent with Armijo line search \citep{Li2018} and with the Polyak step size \citep{You2022} both asymptotically converge for $\alpha \in \interval[open right]{1/2}{1} \cup \interval[open]{1}{\infty}$. 
However, when $\alpha \in \interval[open]{0}{1/2}$, the optimization problem \eqref{eq:renyi_information_measure} becomes non-convex \citep{Lapidoth2019}, and currently, there are no known algorithms that provably solve this problem.
Similarly to the computation of the Augustin information, none of the existing works have established a linear convergence rate.

\section{Preliminaries}\label{sec:preliminaries}

Our analyses are based on properties of Hilbert's projective metric and Birkhoff's contraction theorem, which we introduce in this section.

Let $K$ be a closed cone in a finite-dimensional real vector space, such as the positive orthant and the set of Hermitian positive semidefinite matrices.
For any $x,y\in K$, we write $x\leq y$ if $y-x\in K$.
For any $x,y\in K \setminus \{0\}$, define
\begin{equation}
	M(x/y) \coloneqq \inf\{ \beta\geq 0 \mid x \leq \beta y \} > 0 . 
\end{equation}
If the set is empty, then $M ( x / y ) \coloneqq \infty$.

\begin{definition}
Hilbert's projective metric is defined as
\[
	\dH( x, y ) \coloneqq \log ( M(x/y)M(y/x) ) \in [0, \infty] , \quad \forall x, y \in K \setminus \set{ 0 }.
\]
In addition, $\dH(0,0)$ is defined to be $0$.
\end{definition}


The following lemma shows that $\dH$ is indeed a metric on the set of rays. 

\begin{lemma}\label{lem:projection}
The following properties hold.
\begin{enumerate}[label=(\roman*)]
\item For any $x,y\in K$ and any $\alpha,\beta>0$, we have $\dH( \alpha x, \beta y ) = \dH( x, y )$.
\item We have $\dH(x,y)=0$ if and only if $x=ry$ for some $r>0$.
\end{enumerate}
\end{lemma}

In the rest of the paper, we will only consider the cone $K=\R_+^d$.

\begin{lemma}\label{lem:hilberts_projective_metric}
Consider Hilbert's projective metric $\dH$ on the cone $K=\R_+^d$.
\begin{enumerate}[label=(\roman*)]
\item For any $x,y\in\R_{+}^d\setminus\{0\}$, we have
\begin{equation}
	M(x/y) = \max_{i\in[d]}\frac{ x(i) }{ y(i) },\quad
	\dH(x,y) = \log\max_{i,j\in[d]}\frac{ x(i)y(j) }{ y(i)x(j) }.
\end{equation}
\item 
$(\ri\Delta([d]), \dH)$ is a metric space \citep[Proposition~2.1.1]{Lemmens2012}.
\end{enumerate}
\end{lemma}

Given the second item above, we will measure the errors of both Augustin's fixed-point iteration and \citeauthor{Kamatsuka2024}'s algorithm in terms of Hilbert's projective metric between their iterates and the minimizer. 
The following lemma lists several properties of Hilbert's projective metric, which are direct consequences of Corollary~2.1.4 and Corollary~2.1.5 of \citet{Lemmens2012}.

\begin{lemma}\label{lem:properties}
The following properties hold.
\begin{enumerate}[label=(\roman*)]
	\item $\dH( x^{ r}, y^{ r} ) \leq \abs{r} \dH( x, y )$ for any $x,y\in\R_{+}^d$ and any $r\in\R\setminus\{0\}$.
	\item $\dH( v\odot x, v\odot y ) \leq \dH( x, y )$ for any $x,y,v\in\R_{+}^d$.
	\item $\dH( A x, A y ) \leq \dH( x, y )$ for any $x,y\in\R_{+}^d$ and $A\in\R_+^{d'\times d}$.
\end{enumerate}
\end{lemma}

When the matrix in Lemma~\ref{lem:properties} (iii) is 
entry-wise strictly positive, \citet{Birkhoff1957} showed that linear transformation defined by it 
is 
a contraction. 


\begin{theorem}[\citet{Birkhoff1957}]\label{thm:birkhoff}
Let $A\in\R_{++}^{m\times n}$.
It holds that
\[
	\dH(Ax, Ay) \leq \lambda(A) \cdot \dH(x, y),\quad\forall x,y\in\R_{+}^n,
\]
where
\[
	\lambda(A) \coloneqq \tanh\left(\frac{\delta(A)}{4}\right) < 1,
\]
and
\[
	\delta(A) \coloneqq \log\max_{(i,j),(i',j')\in[m]\times[n]}\frac{ A(i,j)A(i',j') }{ A(i',j) A(i,j') } \geq 0.
\]
\end{theorem}


\section{Linear Rate of Augustin's Fixed-Point Iteration 
}\label{sec:augustin}


In this section, we show that Augustin's fixed-point iteration converges linearly with respect to Hilbert's projective metric for computing the Augustin information of order $\alpha\in(1/2,1)\cup(1,3/2)$.



\subsection{Augustin's Fixed-Point Iteration}

Define the following operators:
    \begin{equation}\label{eq:operator_augustin}
        T_\alpha(x):=\E_{p\sim P} \left[ T_{\alpha,p}(x) \right],
        \quad 
        T_{\alpha,p}(x) := \frac{ p^{\alpha} \odot x^{1-\alpha}}
        {\braket{ p^{\alpha} , x^{1-\alpha} }},\quad\forall p\in\Delta([d]).
    \end{equation}
Augustin's fixed-point iteration \eqref{eq:cover} can be equivalent written as follows:
\begin{enumerate}
    \item Initialize $x_1\in \Delta([d])$.
    \item For all $t\in\N$, compute $x_{t+1}=T_\alpha(x_t)$.
\end{enumerate}

\begin{lemma}[{\citep[Lemma 13]{Nakiboglu2019}}]\label{lem:augustin_optimality_condition}
For $\alpha\in(0, 1)\cup(1,\infty)$, the optimization problem \eqref{eq:augustin_information} has a unique minimizer $x^\star$, which satisfies 
the fixed-point equation 
$x^\star=T_\alpha(x^\star)$.
\end{lemma}

\subsection{Linear Rate Guarantee}

The main result of this section is the following theorem, which bounds the Lipschitz constant of the mapping $T_\alpha$ with respect to Hilbert's projective metric. 
Its proof is 
postponed to 
the next subsection.

\begin{theorem}\label{thm:augustin_contraction}
For $\alpha\in[0,1)\cup(1,\infty)$, we have
\begin{equation}
	\dH(T_\alpha(x), T_\alpha(y))\leq \gamma\cdot \dH(x, y),\quad\forall x,y\in\Delta([d]),
\end{equation}
where $\gamma:=2|1-\alpha|$.
\end{theorem}


Linear convergence of Augustin's fixed-point iteration for $\alpha \in \interval[open]{1/2}{1} \cup \interval[open]{1}{\infty}$ immediately follows. 

\begin{corollary}\label{cor:augustin_linear_rate}
    For any $\alpha\in(1/2, 1)\cup(1,3/2)$, let $x^\star$ be the minimizer of the optimization problem \eqref{eq:augustin_information} and $\{x_t\}$ be the iterates of Augustin's fixed-point iteration \eqref{eq:cover}.
    We have
    \[
        \dH(x_{t+1}, x^\star) \leq \gamma^t \cdot  \dH(x_1, x^\star),
    \]
    for all $t\in\N$, where $\gamma<1$ is defined in Theorem~\ref{thm:augustin_contraction}.
\end{corollary}

\begin{remark}
    Corollary~\ref{cor:augustin_linear_rate} is meaningful only when $\dH(x_1, x^\star)<\infty$.
    Lemma 13 of \citet{Nakiboglu2019} shows that $\E_p[p]\in\ri\Delta([d])$ implies $x^\star\in\ri\Delta_n$.
    In this case, it suffices to choose $x_1\in\ri\Delta([d])$ to ensure
    that 
    $\dH(x_1, x^\star)<\infty$.
\end{remark}

\subsection{Proof of Theorem~\ref{thm:augustin_contraction}}

The proof primarily consists of two steps, which are reflected by the following two lemmas. 
First, we show that the operators $T_{\alpha, p}$ are Lipschitz with respect to Hilbert's projective metric and bound the Lipschitz constant. 
Then, given that $T_\alpha ( \cdot ) = \mathbb{E}_p \left[ T_{\alpha, p} ( \cdot ) \right]$, we prove a general lemma that bounds Hilbert's projective metric between two random probability vectors in terms of Hilbert's projective metric between their realizations, which is of independent interest. 
The proofs of both lemmas are deferred to Appendix~\ref{sec:omitted_proofs}.

%


\begin{lemma}\label{lem:single_contraction}
For any $\alpha\in[0,1)\cup(1,\infty)$ and $p\in\Delta([d])$, 
\[
	\dH(T_{\alpha,p}(x), T_{\alpha,p}(y))
	\leq \abs{1-\alpha} \dH( x, y ),
	\quad\forall x,y\in\Delta([d]).
\]
\end{lemma}

\begin{lemma}\label{lem:distance_of_expectation}
Let $X, Y: \Omega\to\Delta([d])$ be two random probability vectors, 
where $\Omega$ denotes the sample space.
We have
\[
	\dH\left( \E[ X ], \E[ Y ] \right)
	\leq 2 \sup_{\omega\in\Omega} \dH( X(\omega), Y(\omega) ).
\]
\end{lemma}

Theorem \ref{thm:augustin_contraction} follows immediately: 
By Lemma~\ref{lem:distance_of_expectation}, we write 
\[
	\dH( T_\alpha(x), T_\alpha(y) )
	\leq 2\sup_{p\in\Delta([d])} \dH( T_{\alpha,p}(x), T_{\alpha,p}(y) ).
\]
Then, by Lemma~\ref{lem:single_contraction}, we obtain
\[
	\dH( T_\alpha(x), T_\alpha(y) )
	\leq 2\sup_{p\in\Delta([d])} \abs{1-\alpha} \dH( x, y )
	= 2\abs{1-\alpha} \dH( x, y ).
\]
This completes the proof.


\section{Linear Rate of \citeauthor{Kamatsuka2024}'s Algorithm}

In this section, we show that \citeauthor{Kamatsuka2024}'s algorithm converges linearly with respect to Hilbert’s projective metric for computing the R\'enyi information measure of order $\alpha\in[1/2,1)\cup(1,\infty)$.




For convenience, we will view any $p \in \Delta ( [m] \times [n] )$ as a matrix in $\mathbb{R}_+^{m \times n}$ whose entries sum to $1$. 
We will denote Hilbert's projective metric on both $\mathbb{R}_{++}^m$ and $\mathbb{R}_{++}^n$ by $\dH$.
The associated cone should be clear from the context.

\subsection{\citeauthor{Kamatsuka2024}'s Algorithm}

Define the following two operators:
\begin{equation}\label{eq:operator_renyi}
\begin{split}
	U_\alpha(y) &:= \frac { \left( p^{\alpha} y^{1-\alpha} \right)^{1/\alpha} }{ \norm{\left( p^{\alpha} y^{1-\alpha} \right)^{1/\alpha}}_1 },
	\\
	V_\alpha(x) &:= \frac{ \big( (p^{\alpha})^{\top} x^{1-\alpha} \big)^{1/\alpha} }{ \norm{ \big( (p^{\alpha})^{\top} x^{1-\alpha} \big)^{1/\alpha} }_1 }.
\end{split}
\end{equation}
\citeauthor{Kamatsuka2024}'s algorithm \eqref{eq:alternating_minimization} can be equivalently written as follows:
\begin{enumerate}
\item Initialize $x_1\in\Delta([m])$, and compute $y_1 = V_\alpha( x_1 )$.
\item For all $t\in\N$, compute $x_{t+1} = U_\alpha( y_t )$ and $y_{t+1} = V_\alpha( x_{t+1} )$.
\end{enumerate}
This algorithm is inspired by the following lemma \citep[Lemma~16]{Lapidoth2019}.

\begin{lemma}\label{lem:renyi_optimality_condition}
For $\alpha\in[1/2,1)\cup(1,\infty)$, every minimizer $(x^\star,y^\star)$ of the optimization problem \eqref{eq:renyi_information_measure} satisfies $x^\star = U_\alpha(y^\star)$ and $y^\star = V_\alpha(x^\star)$.
\end{lemma}

\subsection{Linear Rate Guarantee}

The following theorem presents a key observation, showing that the operators $U_\alpha$ and $V_\alpha$ have a Lipschitz constant of $\abs{1-\alpha^{-1}}$ with respect to Hilbert's projective metric.
Its proof is postponed to the next subsection.

\begin{theorem}\label{thm:renyi_contraction}
For $\alpha\in[1/2,1)\cup(1,\infty)$, we have
\begin{equation}\label{eq:renyi_contraction}
\begin{alignedat}{2}
	\dH( V_\alpha(x), V_\alpha(x') ) &\leq \gamma'\cdot \dH( x, x' ),
	\quad &&\forall x,x'\in\Delta([m]),\\
	\dH( U_\alpha(y), U_\alpha(y') ) &\leq \gamma'\cdot \dH( y, y' ),
	\quad &&\forall y,y'\in\Delta([n]),
\end{alignedat}
\end{equation}
where 
$\gamma' := \abs{ 1 - ( 1 / \alpha ) }$.
Moreover, if $p\in\R_{++}^{m\times n}$, then the Lipschitz constant $\gamma'$ can be improved to
\[
	\gamma'' := \abs{ 1 - \frac{1}{\alpha} }\lambda( p^\alpha ) < \gamma',
\]
where $\lambda(\cdot)$ is defined in Theorem~\ref{thm:birkhoff}.

\end{theorem}

Theorem~\ref{thm:renyi_contraction} implies the following corollary, showing that the iterative algorithm converges linearly.
The proof of the corollary is deferred to Appendix~\ref{sec:omitted_proofs}.

\begin{corollary}\label{cor:renyi_linear_rate}
Let $(x^\star, y^\star)$ be a minimizer of the optimization problem \eqref{eq:renyi_information_measure} and $\{ (x_t, y_t) \}$ be the iterates of the iterative algorithm \eqref{eq:alternating_minimization}.

\begin{enumerate}[label=(\roman*)]
\item If $\alpha\in(1/2,1)\cup(1,\infty)$, then
\[
\begin{split}
	\dH( x_{t+1}, x^\star )
	&\leq (\gamma')^{2t} \cdot \dH( x_1, x^\star ),\\
	\dH( y_{t+1}, y^\star )
	&\leq (\gamma')^{2t+1} \cdot \dH( x_1, x^\star ),
\end{split}
\]
for all $t\in\N$, where $\gamma'<1$ is defined in Theorem~\ref{thm:renyi_contraction}.

\item If $\alpha\in[1/2,1)\cup(1,\infty)$ and $p\in\R_{++}^{m\times n}$, then
\[
\begin{split}
	\dH( x_{t+1}, x^\star )
	&\leq (\gamma'')^{2t} \cdot \dH( x_1, x^\star ),\\
	\dH( y_{t+1}, y^\star )
	&\leq (\gamma'')^{2t+1} \cdot \dH( x_1, x^\star ),
\end{split}
\]
for all $t\in\N$, where $\gamma''<1$ is defined in Theorem~\ref{thm:renyi_contraction}.
\end{enumerate}
\end{corollary}

\begin{remark}
Corollary~\ref{cor:renyi_linear_rate} is meaningful only when $\dH(x_1, x^\star) < \infty$.
For $\alpha\in[1/2,1)\cup(1,\infty)$, Lemma~\ref{lem:positivity_of_minimizer} in Appendix~\ref{sec:omitted_proofs} ensures that  $x^\star\in\ri\Delta([m])$ whenever $p\in\R_{++}^{m\times n}$ .
In this case, it suffices to choose $x_1\in\ri\Delta([m])$ to ensure $\dH(x_1, x^\star) <\infty$.
\end{remark}

\subsection{Proof of Theorem~\ref{thm:renyi_contraction}}

By Lemma~\ref{lem:projection} and Lemma~\ref{lem:properties} (i), we have
\[
\begin{split}
	\dH( V_\alpha(x), V_\alpha(x') )
	&= \dH\left( ( (p^{\alpha})^\top x^{1-\alpha} )^{1/\alpha},
	 ( (p^{\alpha})^\top (x')^{1-\alpha} )^{1/\alpha} \right) \\
	 &\leq \alpha^{-1} \dH\left( (p^{\alpha})^\top x^{1-\alpha},
	  (p^{\alpha})^\top (x')^{1-\alpha} \right).
\end{split}
\]
By Lemma~\ref{lem:properties} (iii) and (i), we have
\[
	\dH( V_\alpha(x), V_\alpha(x') ) \leq \abs{1-\alpha^{-1}}\dH(x, x').
\]
This proves the first inequality in \eqref{eq:renyi_contraction}.
The second inequality follows from a similar argument.

Assume $p\in\R_{++}^{m\times n}$.
We can apply Birkhoff's contraction theorem (Theorem~\ref{thm:birkhoff}) instead of Lemma~\ref{lem:properties} (iii) to obtain
\[
	\dH( V_\alpha(x), V_\alpha(x') )
	\leq \abs{1-\alpha^{-1}}\lambda\left( (p^{\alpha})^\top \right)\cdot
	\dH(x, x').
\]
The theorem follows by noticing that $\lambda\left( (p^{\alpha})^\top \right) = \lambda(p^{\alpha}) < 1$.


\section{Discussions}

We have proved that Augustin's fixed-point iteration converges at a linear rate for computing the Augustin information of order $\alpha  \in \interval[open]{1/2}{1} \cup \interval[open]{1}{3/2}$, and that \citeauthor{Kamatsuka2024}'s algorithm converges at a linear rate for computing the R\'{e}nyi information measure of order $\alpha \in \interval[open right]{1/2}{1} \cup \interval[open]{1}{\infty}$. 
In contrast, existing results are asymptotic and apply to a narrower range of $\alpha$. 
Our proofs are simple, demonstrating the effectiveness of selecting an appropriate mathematical structure. 

Preliminary numerical experiments indicate that Augustin's fixed-point iteration may converge linearly for $\alpha \in \interval[open]{0}{1} \cup \interval[open]{1}{2}$. 
This observed range is broader than that we have established.
It is natural to explore extending the range of $\alpha$ that admits linear convergence
%

%
%

\acks{
We thank Marco Tomamichel and Rubboli Roberto for discussions.
C.-E.~Tsai, G.-R.~Wang, and Y.-H.~Li are supported by the Young Scholar Fellowship (Einstein Program) of the National Science and Technology Council of Taiwan under grant number NSTC 112-2636-E-002-003, by the 2030 Cross-Generation Young Scholars Program (Excellent Young Scholars) of the National Science and Technology Council of Taiwan under grant number NSTC 112-2628-E-002-019-MY3, by the research project ``Pioneering Research in Forefront Quantum Computing, Learning and Engineering'' of National Taiwan University under grant numbers NTU-CC-112L893406 and NTU-CC-113L891606, and by the Academic Research-Career Development Project (Laurel Research Project) of National Taiwan University under grant numbers NTU-CDP-112L7786 and NTU-CDP-113L7763.

H.-C.~Cheng is supported by the Young Scholar Fellowship (Einstein Program) of the National Science and Technology Council, Taiwan (R.O.C.) under Grants No.~NSTC 112-2636-E-002-009, No.~NSTC 113-2119-M-007-006, No.~NSTC 113-2119-M-001-006, No.~NSTC 113-2124-M-002-003, and No.~NSTC 113-2628-E-002-029 by the Yushan Young Scholar Program of the Ministry of Education, Taiwan (R.O.C.) under Grants No.~NTU-112V1904-4 and by the research project ``Pioneering Research in Forefront Quantum Computing, Learning and Engineering'' of National Taiwan University under Grant No.~NTU-CC-112L893405 and NTU-CC-113L891605. H.-C.~Cheng acknowledges the support from the ``Center for Advanced Computing and Imaging in Biomedicine (NTU-113L900702)'' through The Featured Areas Research Center Program within the framework of the Higher Education Sprout Project by the Ministry of Education (MOE) in Taiwan.}

\bibliography{./refs}

\appendix

\section{Omitted Proofs}\label{sec:omitted_proofs}

\subsection{Proof of Lemma~\ref{lem:single_contraction}}

By Lemma~\ref{lem:projection}, we have
\[
	\dH(T_{\alpha,p}(x), T_{\alpha,p}(y))
	= \dH( p^{\alpha} \odot x^{1-\alpha} , p^{\alpha} \odot y^{1-\alpha} ).
\]
Applying Lemma~\ref{lem:properties} (i) and (ii) gives
\[
	\dH(T_{\alpha,p}(x), T_{\alpha,p}(y))
	\leq \dH( x^{1-\alpha} , y^{1-\alpha} )
	\leq \abs{1-\alpha} \dH( x, y ).
\]
This proves the lemma.

\subsection{Proof of Lemma~\ref{lem:distance_of_expectation}}

We will use the following lemma, whose proof is postponed to the next subsection.

\begin{lemma}\label{lem:expectation}
	Let $X, Y: \Omega\to\Delta([d])$ be two random probability vectors.
	We have
    \[
		M(\E[X] / \E[Y]) \leq \sup_{\omega\in\Omega} M(X(\omega) / Y(\omega) ).
	\]
\end{lemma}



By Lemma~\ref{lem:expectation}, we have
\[
\begin{split}
	M( \E[X] / \E[Y] )
	&\leq \sup_{\omega \in \Omega} M( X(\omega) / Y(\omega) ), \\
	M( \E[Y] / \E[X] )
	&\leq \sup_{\omega \in \Omega} M( Y(\omega) / X(\omega) ).
\end{split}
\]
Then,
\[
	\dH( \E[X] , \E[Y] )
	\leq \sup_{\omega \in \Omega} \log M( X(\omega) / Y(\omega) )
	+ \sup_{\omega \in \Omega} \log M( Y(\omega) / X(\omega) ).
\]
Since $X(\omega), Y(\omega) \in \Delta([d])$, we have
\[
	M( X(\omega) / Y(\omega) ) \geq 1\quad\text{and}\quad
	 M( Y(\omega) / X(\omega) ) \geq 1.
\]
This implies that
\[
\begin{split}
	M( X(\omega) / Y(\omega) )
	&\leq M( X(\omega) / Y(\omega) )\cdot M( Y(\omega) / X(\omega) ) ,\\
	M( Y(\omega) / X(\omega) )
	&\leq M( X(\omega) / Y(\omega) )\cdot M( Y(\omega) / X(\omega) ),
\end{split}
\]
and hence
\[
\begin{split}
	\dH( \E[X] , \E[Y] )
	&\leq \sup_{\omega \in \Omega} \dH( X(\omega), Y(\omega) ) + \sup_{\omega \in \Omega} \dH( X(\omega), Y(\omega) ) \\
	&= 2\sup_{\omega \in \Omega} \dH( X(\omega), Y(\omega) ),
\end{split}
\]
which completes the proof.

\subsection{Proof of Lemma~\ref{lem:expectation}}
    Let $\overline{M}:=\sup_{\omega \in \Omega} M(X(\omega) / Y(\omega) )$.
    We have
	\begin{equation*}
		\overline{M}\E[Y]=\E [\overline{M} Y] \geq \E [M(X/Y) Y] \geq \E [X].
	\end{equation*}
	The lemma follows from the definition of $M(\E[X]/\E[Y])$.

\subsection{Proof of Corollary~\ref{cor:renyi_linear_rate}}

For both (i) and (ii), by Lemma~\ref{lem:renyi_optimality_condition} and Theorem~\ref{thm:renyi_contraction}, we have 
\[
	\dH(x_{t+1}, x^\star)
	= \dH(U_\alpha(y_t), U_\alpha(y^\star))
	\leq \tilde{\gamma} \cdot \dH(y_t, y^\star),
\]
and
\[
	\dH(y_{t+1}, y^\star)
	= \dH(V_\alpha(x_{t+1}), V_\alpha(x^\star))
	\leq \tilde{\gamma} \cdot \dH(x_{t+1}, x^\star),
\]
where $\tilde{\gamma}=\gamma'$ for (i) and $\tilde{\gamma}=\gamma''$ for (ii).
The corollary follows by applying the above two inequalities alternatively.

\subsection{Lemma~\ref{lem:positivity_of_minimizer}}

We prove the following lemma.

\begin{lemma}\label{lem:positivity_of_minimizer}
For $\alpha\in[0,1)\cup(1,\infty)$, let $(x^\star, y^\star)$ be a minimizer of the optimization problem \eqref{eq:renyi_information_measure} and assume $p\in\R_{++}^{m\times n}$.
Then, we have $x\in\ri\Delta([m])$ and $y\in\ri\Delta([n])$.
\end{lemma}

By Lemma~\ref{lem:renyi_optimality_condition}, we have
\[
	x^\star = U_\alpha(y^\star)
	= \frac { \left( p^{\alpha} (y^\star)^{1-\alpha} \right)^{1/\alpha} }{ \norm{\left( p^{\alpha} (y^\star)^{1-\alpha} \right)^{1/\alpha}}_1 }.
\]
Since $p\in\R_{++}^{m\times n}$ and $y^\star\in\Delta([n])$, the vector $p^{\alpha} (y^\star)^{1-\alpha}$ is entry-wise strictly positive.
This implies that $x^\star$ is entry-wise strictly positive, and hence $x^\star\in\ri\Delta([m])$.
To show $y^\star\in\ri\Delta([n])$, consider the equation $y^\star=V_\alpha(x^\star)$ and apply the same argument.
This completes the proof.

\end{document}